# Sur le rang de $J_0(q)$

E. KOWALSKI, Rutgers University,    P. MICHEL, Université d'Orsay

July 31, 1997

**Résumé**

In this paper, we prove an unconditionnal bound for the analytic rank (i.e the order of vanishing at the critical point of the $L$ function) of the new part $J_0^n(q)$, of the jacobian of the modular curve $X_0(q)$. Our main result is the following upper bound: for $q$ prime, one has

$$rank_a(J_0^n(q)) \ll \dim J_0^n(q)$$

where the implied constant is absolute. All previously known non trivials bounds of $rank_a(J_0^n(q))$ assumed the generalized Riemann hypothesis; here, our proof is unconditionnal, and is based firstly on the construction by Perelli and Pomykala of a new test function in the context of Riemann-Weil explicit formulas, and secondly on a density theorem for the zeros of $L$ functions attached to new forms.

## 1 Introduction

Depuis l'article fondamental de Mestre [Me] sur les formules explicitites de Riemann-Weil, l'étude du rang des variétés abéliennes via les propriétés analytiques de leurs fonctions $L$ a suscité de nombreux travaux (voir [Br1, Br2, F-P, Ma, M1, M2, Rm]). Un point commun de ceux-ci est de supposer vraies trois conjectures standard relatives à la fonction $L$ associée à $\mathcal{A}$ convenablement normalisée: on supposait en effet que

1. $L(\mathcal{A}, s)$ admet un prolongement holomorphe à $\mathbf{C}$ avec une équation fonctionnelle qui relie $L(\mathcal{A}, s)$ à $L(1-s, \mathcal{A})$ (depuis les travaux de Wiles et Diamond, on sait que cette hypothèse est vérifiée si $\mathcal{A}$ est un courbe elliptique qui a réduction semi-stable en 3 et 5).

2. L'ordre d'annulation de $L(\mathcal{A}, \int)$ en $1/2$, que l'on note dans la suite $rang_a(\mathcal{A})$ (pour rang analytique) majore $rang_{\mathbf{Q}}\mathcal{A}$ (c'est la moitié de la conjecture de Birch–Swinnerton-Dyer); on dispose d'éléments en faveur de cette conjecture pour les courbes elliptiques modulaires (Cf. les travaux de Gross-Zagier, Kolyvagin, Coates-Wiles, Rubin, Kato).

3. Enfin $L(\mathcal{A}, s)$ vérifie l'hypothèse de Riemann généralisée (les zéros non triviaux de $L(\mathcal{A}, s)$ sont tous situés sur la droite critique ($\Re e\, s = 1/2$)).

Notons que via les travaux de Grothendieck et Deligne les analogues de ces trois hypothèses sont démontrés quand le corps $\mathbf{Q}$ est remplacé par un corps de fonctions sur un corps fini [Br1, M2].

Dans le cas rationnel, la troisième hypothèse est sans doute la plus sérieuse, et Brumer [Br1] fut le premier à suggérer que certains théorèmes de densités de zéros de fonctions $L$ (provenant d'inégalités de "Grand Crible") devraient pouvoir améliorer la situation et permettre de se passer de cette dernière. Cette intuition s'est avérée correcte mais il fallu attendre quelques années et les travaux de Perelli et Pomykala pour la voir se concrétiser: introduisant une toute nouvelle fonction test (que nous appellerons fonction de Perelli-Pomykala) dans les formules explicites de Riemann-Weil ils démontrent le théorème inconditionnel suivant [P-P]:

**Théorème 1.1** *Soit E une courbe elliptique modulaire alors on a l'égalité*

$$\sum_{d \leq D} \mu^2(d) \sum_{\chi_d} rang_a(E_{\chi_d}) = o(D \log D)$$

*où d parcourt les entiers sans facteurs carrés $\leq D$, $\chi_d$ parcourt l'ensemble des caractères primitifs réels de module d et $E_{\chi_d}$ est la courbe elliptique tordue de E par le caractère $\chi_d$.*

Signalons que la borne de Mestre donnerait inconditionnellement $O(D \log D)$ et $O(D)$ en admettant l'hypothèse de Riemann généralisée. Les résultats de [P-P] reposent sur un théorème de densité pour les fonctions $L(E_{\chi_d}, s)$, lui-même provenant d'estimations profondes de Heath-Brown portant sur les sommes associées aux caractères réels (inégalités qui ne sont pas sans rappeler le grand crible classique).

Nous inspirant de leur travail, nous considérons le problème de borner le rang d'une variété abélienne très particulière, à savoir $J_0^n(q)$, la partie nouvelle de la jacobienne $J_0(q)$ de la courbe modulaire $X_0(q)$. Cette variété abélienne définie sur $\mathbf{Q}$ et a pour conducteur $q$. Notre résultat principal est le suivant

**Théorème 1.2** *Supposons q premier, on a la majoration*

(1) $$rang_a(J_0^n(q)) \ll \dim J_0^n(q).$$

Dans le cas général nous obtenons aussi la borne plus faible (mais non triviale):

(2) $$rang_a(J_0^n(q)) \ll \dim J_0^n(q) \log \log q.$$

Mentionnons que la majoration triviale serait en $O(\dim J_0^n(q) \log q)$ et que, en admettant l'hypothèse de Riemann généralisée, Brumer [Br2] a donné la majoration précise (valable pour tout $q$):

(3) $$rang_a(J_0^n(q)) \leq (3/2 + o(1)) \dim J_0^n(q);$$

(nous profitons de cette occasion pour signaler une erreur dans la preuve de la majoration (3) donnée dans [Rm]: il semble en effet que dans la section 6. de *loc. cit.* un terme en $O(N^{-1/2}Tx)$ se soit transformé en passant de la page 274 à la page 275 en un $O(N^{-1/2}Tx^{1/2})$).

Dans l'autre direction, considérant le signe de l'équation fonctionnelle des formes modulaires primitives, il est facile de montrer la minoration suivante valable pour tout $\epsilon > 0$ et $q$ assez grand:

$$rang_a(J_0^n(q)) \geq (\frac{1}{2} - \epsilon) \dim J_0^n(q).$$

D'autre part dans [Rm], Ram Murty a montré la même minoration pour le rang géométrique

$$rang_{\mathbf{Q}}(J_0^n(q)) \geq 1/2 \dim J_0^n(q).$$

Le Théorème 1.2 est bien en fait un théorème en moyenne. En effet, notons $S_2(q)^+$ l'ensemble des formes modulaires primitives ("newforms") de poids 2 et de niveau $q$. Alors, par Eichler-Shimura, on a la factorisation

$$L(J_0^n(q), s) = \prod_{f \in S_2(q)^+} L(f, s),$$

qui est l'analogue de la décomposition de la fonction $L$ d'un corps cyclotomique en produit de fonctions $L$ associées aux caractères de Dirichlet. La majoration (1) peut donc s'écrire ($q$ premier)

(4) $$\sum_{f \in S_2(q)^+} \text{rang}_a(f) = O(|S_2(q)^+|),$$

où $rang_a(f)$ désigne l'ordre d'annulation de $L(f, s)$ en $s = 1$.

A l'instar du Théorème 1.1, les majorations (1), (2) et 4 résultent elles aussi de théorèmes de densité pour les zéros des fonctions $L(f, s)$ qui sont montrés dans le compagnon de cet article [KM].

Pour toute forme $f$, notons $(f, f)$ le produit scalaire de Petersson de $f$ avec elle-même; on peut aussi considèrer le procédé de moyenne suivant sur l'ensemble $S_2(q)^+$ que nous appellerons moyenne "harmonique": posons

$$\sum_{f \in S_2(q)^+}^h X_f := \sum_{f \in S_2(q)^+} \frac{X_f}{4\pi(f, f)}.$$

**Remarque.** — Du point de vue de l'arithmétique, le poids $1/4\pi(f, f)$ qu'on attache à $f$ peut paraitre moins naturel que le poids $1/|S_2(q)^+|$; En revanche de telles moyennes sont parfaitement naturelles dans le cadre de l'analyse harmonique sur $X_0(q)$ voir, par exemple le récent travail [M-U] où c'est la moyenne harmonique qui est la plus naturelle pour la géométrie d'Arakelov de $X_0(q)$. On notera l'égalité (valable pour $q$ premier)

(5) $$\sum_{f \in S_2(q)^+}^h 1 = 1 + O(q^{-3/2})$$

de sorte que les poids $1/4\pi(f, f)$ définissent asymptotiquement une mesure de probabilité sur l'ensemble $S_2(q)^+$.

On démontrera également le théorème suivant qui donne une majoration de ce qu'on pourrait appeler le rang harmonique de $J_0^n(q)$:

**Théorème 1.3** *Supposons $q$ premier, on a la majoration*

$$\sum_{f \in S_2(q)^+}^h \text{rang}_a(f) \ll 1;$$

*la constante impliquée est absolue et explicitable.*

Notons encore que la majoration triviale serait en $\log q$, alors que, admettant l'hypothèse de Riemann généralisée, Ram Murty [Rm] a donné la borne:

$$\sum_{f\in S_2(q)^+}^h \text{rang}_a(f) \leq \frac{7}{6} + o(1).$$

Dans la dernière partie de ce travail, nous améliorons (sous G.R.H) la borne de Ram Murty de manière sensible :

**Théorème 1.4** *Soit $q$ un nombre premier. Supposons que pour chaque forme primitive $f \in S_2(q)^+$, sa fonction $L$ satisfait à l'hypothèse de Riemann. Alors, quand $q \to +\infty$, on a la majoration*
$$\sum_{f\in S_2(q)^+}^h \text{rang}_a(f) \leq \frac{23}{22} + o(1).$$

Cet article n'aurait sans doute pas vu le jour sans les encouragements répétés et les nombreuses suggestions de E. Fouvry et de H. Iwaniec, qu'ils en soient remerciés.

Dans la suite, $\epsilon$ désignera un réel positif aussi petit que l'on souhaite dont la définition peut varier d'une ligne à l'autre. De même pour les constantes "grandes" notées généralement $A$ et $B$.

## 2  Normalisations et lemmes préliminaires

On note $S_2(q)$ (resp. $S_2(q)^+$) l'espace des formes modulaires paraboliques de poids 2 et de niveau $q$ (resp. l'ensemble des formes primitives.)

Toute forme parabolique $f$ a un développement de Fourier (à la pointe infinie)

$$f(z) = \sum_{n\geq 1} \lambda_f(n) n^{1/2}\, \text{e}(nz),\ \text{e}(z) := \exp(2i\pi z).$$

Soit $f \in S_2(q)^+$ une forme primitive, normalisée par $\lambda_f(1) = 1$. Ses coefficients de Fouriers $\lambda_f(n)n^{1/2}$ sont réels.

Pour tout nombre premier $p$ ne divisant pas $q$, on a d'après Deligne-Eichler-Shimura, $\lambda_f(p) = \alpha_p + \overline{\alpha}_p$ avec $|\alpha_p| = 1$; on pose alors pour tout $n \geq 1$, $a_{p^n}(f) = \alpha_p^n + \overline{\alpha}_p^n$. Si $p\|q$, $\lambda_f(p) = \pm 1/p^{1/2}$ et on posera $a_{p^n}(f) = \lambda_f(p)^n$; enfin si $p^2|q$, $\lambda_f(p^n) = \alpha_{p^n} = 0$. Soit $L(f,s)$ la fonction $L$ associée à $f$

$$L(f,s) = \sum_{n\geq 1}\frac{\lambda_f(n)}{n^s} = \prod_{p\|q}(1-\frac{\lambda_f(p)}{p^s})^{-1}\prod_{(p,q)=1}(1-\frac{\alpha_p}{p^{s-1/2}})^{-1}(1-\frac{\overline{\alpha}_p}{p^{s-1/2}})^{-1};$$

Cette fonction admet un prolongement analytique à $\mathbf{C}$ et vérifie l'équation fonctionnelle

$$\Lambda(f,s) = \epsilon_f \Lambda(1-s),\ \epsilon_f = \lambda_f(q)q^{1/2} = \pm 1,\ \Lambda(f,s) = \left(\frac{q}{2\pi}\right)^{s/2}\Gamma(s+1/2)L(f,s)$$

## 2.1 Base orthonormée, sommes de Kloosterman et grand crible

.

Soit $\mathcal{F}$ est une base orthonormée de $S_2(q)$, la "formule des traces" de Petersson relie les coefficients de Fourier des éléments de $\mathcal{F}$ à des sommes de sommes de Kloosterman (cf. [DFI]): pour tout $m, n \geq 1$ on a l'égalité

$$(6) \qquad \frac{1}{4\pi} \sum_{f \in \mathcal{F}} \lambda_f(m) \lambda_f(n) = \delta_{m,n} - 2\pi \sum_{c \geq 1} \frac{S(m,n;cq)}{cq} J_1\left(\frac{4\pi\sqrt{mn}}{cq}\right)$$

où $\delta_{m,n}$ est le symbole de Kronecker,

$$S(m,n;c) = \sum_{\substack{x (\text{mod } c) \\ (x,c)=1}} e(mx + n\overline{x}/c)$$

est la somme de Kloosterman et $J_1(x)$ est la fonction de Bessel d'ordre 1; nous utiliserons les majorations suivantes de $S(m,n;c)$ et $J_1(x)$:

$$(7) \qquad S(m,n;c) \leq (m,n,c)^{1/2} c^{1/2} \tau(c), \quad J_1(x) \ll \min(x, x^{-1/2}).$$

On en déduit l'égalité pour tout $m, n$ (cf. [Du])

$$(8) \qquad \frac{1}{4\pi} \sum_{f \in \mathcal{F}} \lambda_f(m) \overline{\lambda_f(n)} = \delta_{m,n} + O(\tau(q)(m,n,q)^{1/2}(mn)^{1/2} q^{-3/2})$$

**Remarque.** — Un cas particulièrement agréable de cette égalité est quand $q$ est premier: en poids 2, $S_2(q)$ n'a pas de formes anciennes, on peut alors prendre

$$\mathcal{F} = \{\frac{f}{(f,f)^{1/2}}\}_{f \in S_2(q)^+}.$$

## 2.2 Théorèmes de densité

.

Nous citons maintenant les théorèmes de densité pour les fonctions $L(f, s)$ que nous utiliserons dans la section suivante, ils sont démontrés dans l'article compagnon de ce travail [KM]: soit $f \in S_2(q)^+$ pour $1/2 < \alpha \leq 1$ et $t_1 < t_2$ on note $N_f(\alpha, t_1, t_2)$ le nombre de zéros de $L(f, s)$ contenus dans le rectangle $[\alpha, 1] \times [t_1, t_2]$, on posera également $N_f(\alpha, T) = N_f(\alpha, -T, T)$. On a d'abord le résultat général suivant

**Théorème 2.1** *Soit $T \geq 1$, alors il existe des constantes absolues positives $B$ et $c$ telles que l'on ait la majoration*

$$\sum_{f \in S_2(q)^+} N_f(\alpha, T) \ll \dim J_0^n(q) q^{-c(\alpha - 1/2)} T^B \log^B q.$$

Quand $q$ est premier nous améliorons très sensiblement ce résultat en se plaçant dans de petits intervallles et en obtenant $A = 1$:

**Théorème 2.2** *Supposons $q$ premier. Soit $T \geq 1$, et $-T \leq t_1 < t_2 \leq T$ vérifiants $t_2 - t_1 \geq 1/\log q$ alors il existe des constantes absolues positives $B$ et $c$, telles qu' on a la majoration*
$$\sum_{f \in S_2(q)^+} N_f(\alpha, t_1, t_2) \ll T^B q^{1-c(\alpha-1/2)}(t_2 - t_1) \log q.$$

Nous montrons également la version harmonique de ce dernier théorème qui est en quelque sorte en amont du précédent:

**Théorème 2.3** *Supposons $q$ premier. Soit $T \geq 1$, et $-T \leq t_1 < t_2 \leq T$ vérifiants $t_2 - t_1 \geq 1/\log q$ alors il existe des constantes absolues positives $B$ et $c$, telles qu' on a la majoration*
$$\sum_{f \in S_2(q)^+}^h N_f(\alpha, t_1, t_2) \ll T^B q^{-c(\alpha-1/2)}(t_2 - t_1) \log q.$$

### 2.3 Les formules explicites de Riemann-Weil

.

Considérons $f \in S_2(q)^+$; la proposition suivante est due à Mestre [Me]:

**Proposition 2.4** *Soit $F : \mathbf{R} \to \mathbf{R}$ vérifiant*

1. *Il existe $\epsilon > 0$ tel que $F(x)\exp((1+\epsilon)x)$ est intégrable et à variation bornée,*

2. *la fonction $(F(x) - F(0))/x$ est à variation bornée.*

*On a la formule suivante,*

$$\sum_{L(f,\rho)=0} \widehat{F}(\rho - 1/2)$$
$$= 2F(0)\log q^{1/2} - 2\sum_{n=1}^{\infty} \frac{a_n(f)}{n^{1/2}}\Lambda(n)F(\log n) - \int_0^{+\infty}\left(\frac{F(x)\mathrm{e}^{-x}}{1 - \mathrm{e}^{-x}} - \frac{\mathrm{e}^{-x}}{x}\right)dx.$$

*On a posé*
$$\widehat{F}(s) = \int_{\mathbf{R}} F(x)\mathrm{e}^{sx}\,dx.$$

*La somme du terme de gauche porte sur les zéros de $L(f,s)$ tels que $0 \leq \Re e\rho \leq 1$.*

**Remarque.** — Pour toutes les fonctions test $F$ que nous utiliserons dans la suite la dernière intégrale de l'expression précédente est bornée par une constante absolue (cf. [P-P] Section 5).

Le fait fondamental qui est à l'origine de cet article est la construction par Perelli et Pomykala d'une fonction test $F_{PP}$ aux propriétés extrêmement agréables quand on doit se passer de l'hypothèse de Riemann généralisée [P-P]:

**Théorème 2.1** *Il existe des constantes positives $B, c_1, c_2$ avec $B \geq 2$, une fonction réelle paire à valeurs non négatives $F_{PP} \in \mathcal{C}^\infty(\mathbf{R})$ telle que $F_{PP}(0) = 1$, $\operatorname{supp} F_{PP} \subset [-B, B]$, telle que $\widehat{F}_{PP}(s)$ vérifie*

$$\Re e(\widehat{F}_{PP}(s)) \geq 0, \text{ pour } |\Re e s| < 1, \text{ et}$$

$$\widehat{F}_{PP}(s) \ll \exp(c_1 |\Re e\ s| - c_2 |s|^{3/4}).$$

## 3 Majorations inconditionnelles du rang analytique de $J_0^n(q)$

La méthode pour majorer la quantité $\sum_{f \in S_2(q)^+} \operatorname{rang}_a(f)$ suit de près celle de Perelli et Pomykala [P-P]. On se contentera ici de montrer la majoration (1), la preuve de la borne (2) est identique, il suffit d'utiliser le Théorème 2.1 à la place du Théorème 2.2

### 3.1 Majoration du Rang analytique

.

Soit $\lambda > 2$, on pose $F_\lambda(x) := F_{PP}(x/\lambda)$ et $\phi_\lambda(s) := \widehat{F}_\lambda(s - 1/2)$ $= \lambda \widehat{F}_{PP}(\lambda(s - 1/2))$; on a en particulier $\phi_\lambda(s) = \phi_\lambda(2 - s)$. Soit $T \geq 2$, un paramètre à fixer. Pour chaque $f$ primitive, on applique la formule de Riemann-Weil avec $F_\lambda$ pour fonction test, les propriétées de $F_{PP}$ et la majoration $N_f(1, t) \ll t \log(qt)$ pour $t \geq T$ permettent d'en déduire l'égalité

$$\sum_{|\gamma| \leq T} \phi_\lambda(\rho) = \log q - 2S_{1,f}(\lambda) - 2S_{2,f}(\lambda) +$$

(9)
$$O(\lambda T \log(qT) \exp(c_5 \lambda - c_6 (\lambda T)^{3/4}))$$

avec

$$S_{1,f} = \sum_p a_p(f) \frac{F_\lambda(\log p) \log p}{p^{1/2}}, \text{ et}$$

$$S_{2,f} = \sum_{p, n \geq 2} a_{p^n}(f) \frac{F_\lambda(n \log p) \log p}{p^{n/2}}.$$

Par la majoration $|a_p^n(f)| \leq 2$ on a immédiatement la majoration

$$S_2(\lambda) = O(\lambda).$$

D'autre part, Brumer ([Br2] Prop 2.8), utilisant la formule des traces de Selberg dans la version qu'en a donné Zagier, obtient la majoration suivante pour tout nombre premier $p$ ($\tau(n)$ est la fonction nombre de diviseurs):

$$\sum_f \lambda_f(p) = O(\tau(q)^3 \log^2(pq)(p^{1/2} + q^{1/2})),$$

dont on déduit la majoration

(10) $$\sum_f S_{1,f} = O(\lambda \tau(q)^3 \log^2(q)(e^{B\lambda} + e^{B\lambda/2} q^{1/2})).$$

Sommant (9) sur les $f$ et prenant la partie réelle, on obtient alors la majoration

$$\lambda \sum_f \mathrm{rang}_a(f) \ll \dim J_0^n(q)(\log q + \lambda) + \lambda \tau(q)^3 \log^2(q)(e^{B\lambda} + e^{B\lambda/2} q^{1/2})$$
$$+ q\lambda T \log(qT) \exp(c_5 \lambda - c_6(\lambda T)^{3/4}) + \sum_f \sum_{\substack{|\gamma| \leq T \\ |\beta - 1/2| \geq 1/\lambda}} |\phi_\lambda(\rho)|$$

Pour majorer la dernière somme de l'expression précédente, on subdivise les intervalles $[1/\lambda, 1/2]$ et $[-T, T]$ en sous intervalles de longueur $1/\lambda$, si $q$ est premier, on a par le Théorème 2.2

$$\sum_f \sum_{\substack{|\gamma| \leq T \\ |\beta - 1/2| \geq 1/\lambda}} |\phi_\lambda(\rho)|$$
$$\ll \lambda \dim J_0^n(q) \sum_{m=1}^{\lambda} \sum_{n=0}^{T\lambda} \exp(c_1 m - c_2 n^{3/4})(\frac{n}{\lambda} + 1)^B q^{-cm/\lambda} \frac{\log q}{\lambda}$$
$$\ll \lambda \dim J_0^n(q),$$

en choisissant $\lambda = c_7 \log q$ ou $c_7$ est une constante positive suffisament petite; On prend $T = \log^3 q$ et on choisit également $c_7$ de sorte que

$$\lambda q T \log(qT) \exp(c_5 \lambda - c_6(\lambda T)^{3/4}) = o(\lambda \dim J_0^n(q)).$$

### 3.2 Majoration du rang harmonique

.

La preuve du Théorème 1.3 suit de près celle de (1): Utilisant l'égalité (5) on obtient la majoration

$$\lambda \sum_f^h \mathrm{rang}_a(f) \ll \log q + \lambda + \sum_f^h S_{1,f} + \lambda T \log(qT) \exp(c_5 \lambda - c_6(\lambda T)^{3/4})$$
$$+ \sum_f^h \sum_{\substack{|\gamma| \leq T \\ |\beta - 1/2| \geq 1/\lambda}} |\phi_\lambda(\rho)|$$

Etablissons l'analogue pondéré de la majoration (10): on a, par (8) pour $m = 1, n = p \leq e^{B\lambda}$,

$$\sum_f^h S_{1,f} \ll \frac{e^{B\lambda}}{q^{3/2}}.$$

utilisant le Théorème 2.3, on majore dernière la somme

$$\sum_{f}^{h} \sum_{\substack{|\gamma| \leq T \\ |\beta - 1/2| \geq 1/\lambda}} |\phi_\lambda(\rho)|$$

$$\ll \quad \lambda \sum_{m=1}^{\lambda} \sum_{n=0}^{T\lambda} \exp(c_1 m - c_2 n^{3/4})(\frac{n}{\lambda} + 1)^B q^{-cm/\lambda} \frac{\log q}{\lambda}$$

$$\ll \quad \lambda.$$

On obtient cette dernière majoration en choisissant $\lambda = c_8 \log q$ ou $c_8$ est une constante positive suffisament petite. On termine alors la preuve comme précédemment.

## 4  Majoration sous G.R.H

Dans cette section on suppose que $q$ est premier et que l'hypothèse de Riemann généralisée est vérifiée pour toutes les fonction $L(f, s)$, $f \in S_2(q)^+$. Grâce à cette lourde hypothèse il est loisible de choisir une fonction test classique [Rm, F-P, M1, M2, Br1]: pour $\lambda > 0$ on pose

$$F_\lambda(x) = \max(0, 1 - |x|/\lambda) \text{ et } \widehat{F}_\lambda(i\gamma) = \phi_\lambda(i\gamma) = \frac{4\sin^2(\gamma\lambda/2)}{\lambda\gamma^2};$$

Sous GRH, $\phi_\lambda(\rho - 1) \geq 0$ pour tout zéro non trivial $\rho$ de $L(f, s)$, on a donc la majoration

(11) $$\lambda \operatorname{rang}_a(f) \leq \log q - 2S_{1,f} - 2S_{2,f} + O(1)$$

On fait la somme sur les formes primitives $f$ pour obtenir grâce à (5)

$$\lambda \sum_{f \in S_2(q)^+}^{h} \operatorname{rang}_a(f) \leq \log q(1 + O(q^{-3/2})) - 2\sum_{f}^{h} S_{1,f} + \sum_{f}^{h} S_{2,f} + O(1).$$

On traite maintenant le troisième terme par une méthode, différente de celle de [Rm], et beaucoup plus élémentaire puisqu'elle ne fait pas intervenir la Théorie de Rankin-Selberg. Notons que par $|a_{p^n}(f)| \leq 2$, on a

(12) $$\sum_{\substack{p \\ n \geq 3}} a_{p^n}(f) \frac{F_\lambda(n \log p) \log p}{p^{n/2}} = O(1).$$

Pour $n = 2$, on part de l'égalité $a_{p^2}(f) = a_p(f)^2 - 2$, puis utilisant (6) puis l'estimation (8) pour $m = n = p$,

$$-2\sum_{f}^{h} S_{2,f} = -2 \sum_{p \leq e^{\lambda/2}} \frac{F(2\log p)\log p}{p}(-1 + O((q,p)^{1/2}pq^{-3/2})) + O(1)$$

(13) $$= \frac{\lambda}{2} + O(e^{\lambda/2 - 3/2 \log q}) + O(1)$$

### 4.1 Traitement du deuxième terme

.

C'est pour ce terme que l'on améliore le Théorème 1 de [Rm]: au lieu d'utiliser une majoration individuelle des sommes de Kloosterman, on tient compte de leur variation quand la variable $p$ décrit l'ensemble des nombres premiers. On a les égalités suivantes

$$\sum_f^h S_{1,f} = \sum_{c \geq 1} \frac{1}{cq} \sum_p \frac{F_\lambda(\log p) \log p}{p^{1/2}} S(1,p;cq) J_1(\frac{4\pi\sqrt{p}}{cq})$$

$$:= \sum_{c \geq 1} \frac{1}{cq} \Sigma_{cq}.$$

Pour $c > C$, on utilise (7) pour majorer $\Sigma_{cq}$ par $O_\epsilon(\lambda(Cq)^{-3/2+\epsilon})$. On suppose dans la suite que $1 \leq c \leq C$. Par intégration par partie et en utilisant la majoration $J_1'(x) = O(1)$, on obtient que

$$\Sigma_{cq} \ll \frac{\lambda}{cq} \max_{1 \leq x \leq X} \Sigma'(x),$$

où on a posé $X = e^\lambda$ et

$$\Sigma'(x) = \sum_{n \leq x} \Lambda(n) S(1,n;cq) + O(x^{1/2}(cq)^{1/2+\epsilon}).$$

Le traitement de cette somme est calqué sur celui de Vaughan [Va] Corollary 2.1:

$$\Sigma'(x) = \sum_{a=0}^{cq-1} S(1,a;cq) \sum_{\substack{n \leq x \\ n \equiv a(cq)}} \Lambda(n)$$

$$= \frac{1}{\phi(cq)} \sum_{\substack{a=1 \\ (a,cq)=1}}^{cq-1} S(1,a;cq) \sum_{\chi (\text{mod } cq)} \overline{\chi}(a) \psi(x,\chi) + O_\epsilon((cq)^{1/2}\epsilon) \log x)$$

$$\ll \frac{1}{\phi(cq)} \sum_{\chi (\text{mod } cq)} \Big| \sum_{\substack{a=1 \\ (a,cq)=1}}^{cq-1} S(1,a;cq) \overline{\chi}(a) \Big| |\psi(x,\chi)| + O((cq)^{1/2+\epsilon}) \log x).$$

On utilise alors le lemme facile suivant:

**Lemme 4.1** *Soit $\chi$ un caractère modulo $cq$, on a la majoration*

$$\sum_{\substack{a=1 \\ (a,cq)=1}}^{cq-1} S(1,a;cq) \chi(a) = O(cq).$$

Ce lemme se démontre en ouvrant la somme de Kloosterman ce qui nous ramène à des sommes de Gauss, puis en intervertissant les sommations...

Ce lemme et le corollaire 2 de [Va] conduisent à la la majoration

$$|\Sigma'(x)| = O_\epsilon(x + (cq)^{5/8} x^{3/4} + cq x^{1/2}).$$

Utilisant cette dernière, on obtient

$$\text{(14)} \qquad \sideset{}{^h}\sum_f S_{1,f} \ll_\epsilon (cq)^\epsilon \left( \frac{X}{q^2} + \frac{X^{3/4}}{q^{11/8}} + \frac{X^{1/2} \log C}{q} + \frac{X}{(Cq)^{3/2}} \right).$$

On choisit alors $C = q^{1/2}$, $\lambda = (11/6 + \epsilon) \log q$ si bien qu'en réunissant (14), (13) et (11) on obtient la majoration

$$\sideset{}{^h}\sum_{f \in S_2(q)^+} \text{rang}_a(f) \leq \frac{6}{11} + \frac{1}{2} + o(1).$$

# Références

[Br1]      A. BRUMER. — *The average rank of elliptic curves I*, Invent.Math. 109, (1992) , 445-472.

[Br2]      A. BRUMER. — *The rank of $J_0(N)$*, Astérisque 228, SMF (1995) , 41-68.

[DFI]      W. DUKE, J.B. FRIEDLANDER et H. IWANIEC. — *Bounds for Automorphic L-functions II*, Invent. Math 115, 219-239 (1994).

[Du]      W. DUKE. — *The critical order of vanishing of automorphic L-functions with high level* Invent. Math. 119, (1995), 165-174.

[F-P]      E. FOUVRY et J. POMYKALA. — *Rang des courbes elliptiques et sommes d'exponentielles,* Mh. Math. 116 (1993), 111-125.

[Iw]      H. IWANIEC. — *Topics in classical automorphic forms*, Grad. Studies in Math., A.M.S, à paraître.

[KM]      E. KOWALSKI et P. MICHEL. — *Sur les zéros des fonctions L automorphes de grand niveau*, préprint.

[Ma]      L. MAI *The average analytic rank of a family of elliptic curves,* J. Number Theory 45, No.1, 45-60 (1993).

[Me]      J.-F. MESTRE. — *Formules explicites et minorations de conducteurs de variétés algébriques,* Comp. Math. 58, (1986), 209-232.

[M1]      P. MICHEL. — *Rang moyen de famille de courbes elliptiques et lois de Sato-Tate.* Mh. Math. 120, (1995) 127-136.

[M2]      P. MICHEL. — *Rang moyen de famille de variétés abéliennes.* J. Alg. Geom., vol. 6 n. 2 (1997), 201-234.

[P-P]      A. PERELLI et J. POMYKALA. — *Averages over twisted elliptic L-functions,* Acta Arith. 80, No 2, 149-163 (1997).

[M-M]      L. MAI et R. MURTY. — *The Phragmen-Lindelöf theorem and modular elliptic curves,* Contemporary Math., 166 (1994) p. 335-340.

[M-U]      P. MICHEL et E. ULLMO. — *Sur les points de petite hauteur des courbes modulaires $X_0(N)$,* Inventiones Mathematicae (à paraitre).

[Rm]      R. MURTY. — *The analytic rang of $J_0(N)$,* (K. Dilcher Edt.) CMS Conf. Proc. 15, 263-277 (1995).

[Va]      R.C. VAUGHAN. — *Mean Value Theorems in Prime Number Theory,* J. London Math. Soc. (2) (10), (1975), 153-162.

P. MICHEL: michel@math.u-psud.fr
Université PARIS-SUD, Bât. 425 91405 ORSAY Cdx. FRANCE

E. KOWALSKI: ekowalsk@math.rutgers.edu
Dept. of Math. RUTGERS University, New Brunswick, NJ 08903 USA